\documentclass[11pt]{amsart}
\usepackage[T1]{fontenc}
\usepackage{latexsym,amssymb,amsmath}
\input epsf

\def\cM{{\mathcal M}}
\def\cN{{\mathcal N}}
\def\cS{{\mathcal S}}
\def\FF{\mathbb F}

\def\lra{\longrightarrow}

\newtheorem*{theo}{Theorem}
\newtheorem*{coro}{Corollary}
\newtheorem*{lemm}{Lemma}

\begin{document}

\title{Generic singularities of Schubert varieties}
\author{L. Manivel}
\date{May 2001}

\medskip
\begin{abstract}
We describe the generic singularity of a Schubert variety of type A
on each irreducible component of its singular locus. This singularity 
is given either by a cone of rank one matrices, or a quadratic cone.
\end{abstract}
\maketitle

\section{Introduction}

Let $\FF_n$ be the variety of complete flags of an $n$-dimensional 
vector space over an algebraically closed field of characteristic
zero. For each reference flag $V_{\bullet}$, and each permutation 
$w\in\cS_n$, one can define a Schubert variety 
$$X_w=\{W_{\bullet}\in \FF_n,
\;\dim (W_p\cap V_q)\ge r_w(p,q),\;1\leq p,q\leq n\},$$
where the rank function $r_w$ is defined 
$r_w(p,q)=\# \{i\leq p,\; w(i)\leq q\}.$ This Schubert variety 
is the Zariski closure of an affine cell $\Omega_w$ of dimension 
$l(w)$, the length of the permutation $w$. 

This note is an appendix to \cite{m2}, where we solved the 
longstanding problem of locating the irreducible components 
of Schubert varieties (see also \cite{bw} and \cite{klr}). 
Here we describe geometrically the 
singularity of $X_w$ at the generic point of each component of its 
sigular locus. This was done in \cite{ac} for {\em covexillary} 
permutations (no $3412$ configuration allowed), and in a more
general but also more restricted context (corresponding in type A 
to the case of grassmannians rather than complete flag varieties),
in \cite{bp}. 

Before stating our result we need to recall how the singular locus
of $X_w$ can be located. Recall that following a theorem of Lakshmibai 
and Sandhya, $X_w$ is singular if and only if there is no sequence of
integers $i<j<k<l$ such that $w(l)<w(j)<w(k)<w(i)$ (a $4231$
configuration) or $w(k)<w(l)<w(i)<w(j)$ (a $3412$ configuration).
In general, the irreducible components of ${\rm Sing}(X_w)$ are in
correspondance with certain minimal $4213$ or $3412$  configurations 
encoded in the following figures, where the $\bullet$ are points of
the diagram of $w$ (their coordinates are $(i, w(i))$ for some $i$).
There are actually three types. 
\medskip

\noindent Type $4231$ : $l(w)-l(v) = l+m+1$, 
$m(w,v)= l(w)+lm=l(v)+(l+1)(m+1)$, where $l>0$ and $m>0$ are the numbers
of $\bullet$ in the NorthWest and SouthEast squares, respectively.

$$
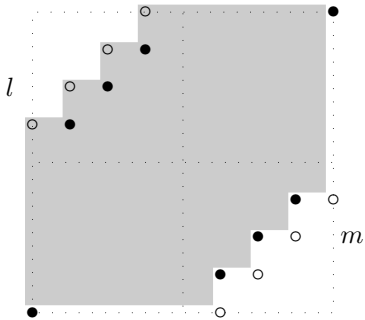
\epsfbox{sing.1}$$

\pagebreak

\noindent Type $34\hspace{-1mm}*\hspace{-1mm}12$ : $l(w)-l(v) = 
2l+3$, $m(w,v)= l(w)+1=l(v)+2l+4$, where $l\ge 0$ is the number of
$\bullet$ in the central square. 

 $$
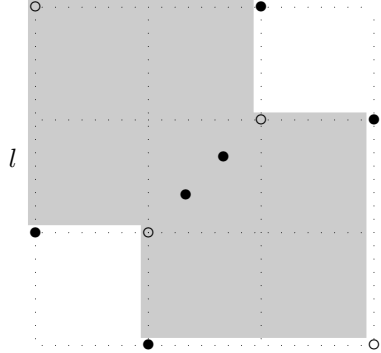
\epsfbox{sing.3}$$

\noindent Type $34\emptyset 12$ : $l(w)-l(v)  = l+m+3$,
$m(w,v)=  l(w)+l+m+1=l(v)+2(l+m+2)$, where $l\ge 0$ and $m\ge 0$ are 
the numbers of $\bullet$ in the NorthWest and SouthEast squares, 
respectively.

$$
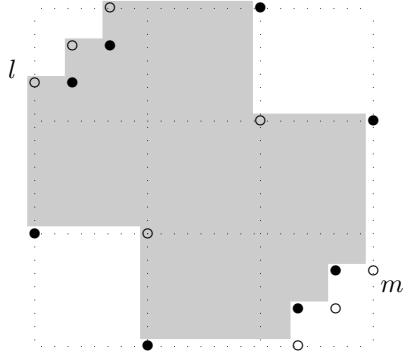
\epsfbox{sing.2}$$

\medskip
The configurations formed by the $\bullet$ in the figures above are
minimal if there is no other point of the diagram of $w$ in  
the region $D$ drawn in grey. For each such configuration, we then
replace the $\bullet$ by the points represented by $\circ$,
to obtain the diagram of a new permutation $v$. Then $X_v$ is 
an irreducible component of the singular locus of $X_w$, and every 
irreducible component is obtained that way. (Moreover, $D$ is 
precisely the region where $r_v>r_w$.) The dimension $m(w,v)$ of 
the Zariski tangent space $T_xX_w$ at a point $x\in\Omega_v$
was computed in \cite{m2}. 
\medskip

\begin{theo}
Let $X_v$ be an irreducible component of $X_w$, coming from one 
of the three possible types of minimal configurations listed above. 
Then each point of $\Omega_v$ has an affine neighbourhood in $\FF_n$ 
whose intersection with $X_w$ is isomorphic to the product of the 
affine cell $\Omega_v$, of dimension $l(v)$, with either
\begin{enumerate}
\item a cone of matrices of size $(l+1)\times (m+1)$ and rank 
at most one;
\item a quadratic cone of dimension $2l+3$;
\item a cone of matrices of size $2\times (l+m+2)$ and rank 
at most one.
\end{enumerate}
\end{theo}

The following corollary is an immediate consequence of the theorem 
and of the computations of \cite{bp}, 3.3. It was obtained in a purely 
combinatorial way in \cite{bw}, 12, but our geometric statement is
of course more precise:

\begin{coro} With the same notations as above, the Kazhdan-Lusztig
  polynomial of the pair $(v,w)$ is, respectively,
\begin{enumerate}
\item $P_{v,w}(q) = 1+q+\cdots +q^{\min(l,m)}$;
\item $P_{v,w}(q) = 1+q^{l+1}$;
\item $P_{v,w}(q) = 1+q$.
\end{enumerate}
\end{coro}

\section{Proof of the theorem}

As in \cite{bp,ac}, we use the existence of a transversal $\cN_{v,w}$ 
to $\Omega_v$ in $X_w$, which was noticed in \cite{kl}, Lemma A.4. 
This transversal (whose dimension is of course $l(w)-l(v)$) is the 
intersection of $X_w$ with $\cN_v=v(\Omega_{w_0})
\cup\Omega_{w_0}$, where $w_0$ is the permutation with maximal length
in $\cS_n$ ($\Omega_{w_0}$ is the ``big cell'' in $\FF_n$, isomorphic
to the unipotent group $U^-$ of strict lower triangular matrices in 
$GL_n$), and the permutation $v$ is identified with its matrix
(assuming that the reference flag is just the canonical flag).
The columns of $v$ generate a flag $D_v\in\Omega_v$, and the 
Bruhat decomposition implies that the map  $\phi_v : \cM_v=vU^-\cap U^-v
\lra \cN_v$ given by $\phi_v(m)=m(D_v)$, is an isomorphism. 
Note that
$$\cM_v = \{m\in GL_n,\; m_{iv(i)}=1, \; m_{jk}=0 \;{\rm if}\; k<v(j)
\;{\rm or} \;j>v^{-1}(k)\}.$$

\smallskip
\begin{lemm} Let $m\in \phi_v^{-1}(\cN_{v,w})$. Let $j,k$ be such that
$k>v(j)$ and $j<v^{-1}(k)$. Then $m_{jk}=0$ as soon as the rectangle 
$[j,v^{-1}(k)[\times [v(j),k[$ is not contained in the region $D$ 
where $r_v>r_w$.\end{lemm}

\proof Choose a basis $e_1,\ldots ,e_n$ adapted to the reference flag 
$V_{\bullet}$. Let $W_{\bullet} =\phi_v(m)$, with $m\in\cM_v$, and 
$(p,q)\notin D$. 
Then $W_p+V_q$ is generated by the vectors $e_1,\ldots ,e_q$ and 
$m(e_1),\ldots ,m(e_p)$. Since $m(e_j)=e_{v(j)}+\sum_{k>v(j)}m_{jk}e_k$,
the subfamily formed by  $e_1,\ldots ,e_q$ and those $m(e_j)$ such
that $j\le p$ and $v(j)>q$ consists in independant vectors. Since
there are $q+p-r_v(p,q)=q+p-r_w(p,q)$ of them, we get $dim(W_p\cap
V_q)\le r_w(p,q)$. If $W_{\bullet}$ belongs to $X_w$, we must have
equality, and this implies that the other vectors of the family, 
that is the $m(e_j)$ for which $j\le p$ and $v(j)\le q$, must be
linear combinations of the previous ones. Hence $m_{jk}=0$ if 
$j\le p$, $v(j)>q$, and $k>q$ is not among $v(1),\ldots ,v(p)$,
that is $v^{-1}(k)>p$. The lemma follows immediately. \qed

\medskip Now we study our three cases separately. 

\medskip\noindent {\em First case}. The first one,
that of a minimal $4231$ configuration, is completely similar to 
\cite{ac}, Th\'eor\`eme 3.6. We fix our notations as in the 
figure below. 

$$
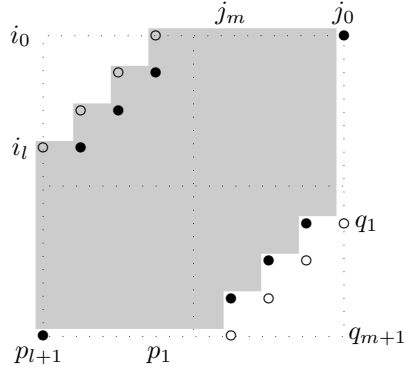
\epsfbox{sing.4}$$

The Lemma implies that a matrix $m\in\phi_v^{-1}(\cN_{v,w})$ can 
have non zero entries, except for those which must be equal to one,
only on the lines $j_0,\ldots ,j_l$ and the columns $k_0,\ldots,
k_m$. Denote by $J$ and $K$ these sets of indices, this makes an 
$(l+1)\times (m+1)$ submatrix $m_{J,K}$ of indeterminates. 

Now consider the incidence condition corresponding to the point 
$(i_l,p_1)$, which belongs to $D$. As in the proof of the Lemma,
we see that $W_{i_l}+V_{p_1}$ contains the independant family
consisting in the vectors $e_1,\ldots ,e_{p_1}$ and the $m(e_j)$ 
for $j\le i_l$ and $v(j)>p_1$. Moreover, since $r_v(i_l,p_1)=
r_w(i_l,p_1)+1$, the dimension of $W_{i_l}+V_{p_1}$ is at most 
one more than the number of vectors in this family. 
But $W_{i_l}+V_{p_1}$ also contains $m(e_{i_0}),\ldots ,m(e_{i_l})$, 
and in consequence, the rank of $M$ cannot be larger than one. 
Therefore $\phi_v^{-1}(\cN_{v,w})\subset \cN'_{v,w}$, where 
$$\cN'_{v,w}=\{m\in GL_n,m_{iv(i)}=1, \; m_{jk}=0 \; {\rm if} j\notin
J\;{\rm or}\; k\notin K, \; {\rm rank}(m_{J,K})\le 1\}.$$
But this is an irreducible variety of the same dimension as 
$\cN_{v,w}$, hence there must be equality. 

\medskip

\pagebreak\noindent {\em Second case}.
We fix our notations as in the figure below. 

$$
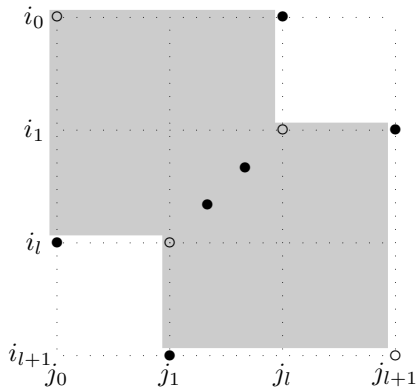
\epsfbox{sing.6}$$
 
Here the Lemma implies that a matrix $m\in\phi_v^{-1}(\cN_{v,w})$ can 
have non zero entries, except for those which must be equal to one,
only on the line $j_0$ and the columns $k_1,\ldots, k_{l+1}$, 
or on the column $k_{l+2}$, and the lines $j_1,\ldots, j_{l+1}$. 

Consider the incidence condition corresponding to the point 
$(j_l,k_0)$, which does not belong to $D$. It implies that 
$W_{j_l}+V_{k_0}$ has a basis consisting of $e_1,\ldots ,e_{k_0}$
and the $m(e_j)$ such that $j\le j_l$ and $v(j)>k_0$. Then
$$m(e_{j_0})=e_{k_0}+a_1e_{k_1}+\cdots +a_{l+1}e_{k_{l+1}}$$ 
must be a linear combination of these vectors, hence 
of $m(e_{j_1})=e_{k_{l+1}}+b_1e_{k_{l+2}},\ldots ,
m(e_{j_{l+1}})=e_{k_1}+b_{l+1}e_{k_{l+2}}$, and $e_{k_0}$.
This is equivalent to the quadratic condition $a_1b_{l+1}+\cdots
+a_{l+1}b_1=0$. 

All these conditions define a quadratic cone $\cN'_{v,w}$
containing $\phi_v^{-1}(\cN_{v,w})$. Since it is irreducible of the 
same dimension as $\cN_{v,w}$, again there must be equality. 

\medskip\noindent {\em Third case}.
This case is slightly more complicated than the previous ones. 
We fix our notations as in the figure below. 

$$
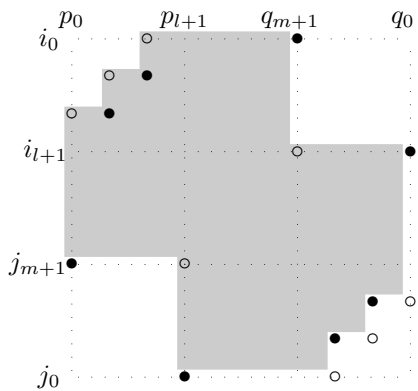
\epsfbox{sing.5}$$
 
The Lemma implies that a matrix $m\in\phi_v^{-1}(\cN_{v,w})$ can 
have non zero entries, except for those which must be equal to one,
only on the lines $i_0\ldots ,i_l$ and the columns $p_{l+1}, q_{m+1}$
(giving the coefficients of a $(l+1)\times 2$ matrix $A$), 
or on the lines $i_{l+1}, j_{m+1}$  and the columns $q_0,\ldots, q_m$
(giving the coefficients of a $2\times (m+1)$ matrix $B$). 
 
Exactly as in the first case, the incidence conditions corresponding 
to the points $(i_l,p_l)$ and $(j_{m+1}, q_{m+1})$  impose that 
${\rm rank}(A)\le 1$ and ${\rm rank}(B)\le 1$, respectively. 

Now consider the incidence condition given by the point $(j_{m+1},
p_l)$, which is not in $D$. It implies that 
$W_{j_{m+1}}+V_{p_l}$ has a basis consisting of $e_1,\ldots ,e_{p_l}$
and the $m(e_j)$ such that $j\le j_{m+1}$ and $v(j)>p_l$. Then
$m(e_{i_0}),\ldots , m(e_{i_l})$ must be linear combinations of the 
previous vectors. Thus, more precisely, $m(e_{i_j})$ must be a linear
combination of $e_{p_{l-j}}$, $m(e_{i_{l+1}})$ and $m(e_{j_{m+1}})$. 
If we denote by $a, a'$ the two columns of $A$, and by $b, b'$ the 
two lines of $B$, this means that $a_ib_j+a'_ib'_j=0$ for all $i,j$. 

All these conditions define an irreducible variety $\cN'_{v,w}$
containing $\phi_v^{-1}(\cN_{v,w})$. Since it has  the 
same dimension as $\cN_{v,w}$, again there must be equality. 

Note that if we see the matrices $A$ and $B$ as morphisms
$A : k^2\lra k^{l+1}$ and $B : k^{m+1}\lra k^2$, they must have 
rank at most one and $A\circ B$ must be zero. If we associate to
such a pair the $2\times (l+m+2)$ matrix
$$ C = \begin{pmatrix} a_0& \cdots & a_{l+1} & b'_0 & \cdots &
  b'_{m+1} \\ a'_0& \cdots & a'_{l+1} & -b_0 & \cdots &
  -b_{m+1} \end{pmatrix},$$
we obtain an isomorphism on the cone of rank one matrices, and the 
proof of the Theorem is complete.

\bigskip\noindent
Laurent Manivel, {\sc Institut Fourier}, UMR 5582 du CNRS, Université 
Joseph Fourier, BP 74, 38402 Saint Martin d'Hères, France. 

\smallskip
\noindent {\em E-mail} : Laurent.Manivel@ujf-grenoble.fr


\bigskip
\begin{thebibliography}{aa}

\bibitem{bw} Billey S. C., Warrington G. S., Maximal singular loci of
  Schubert varieties in $SL(n)/B$, preprint arXiv:math.AG/0102168.

\bibitem{bp} Brion M., Polo P., Generic singularities of certain 
Schubert varieties, {\sl Math. Z.} {\bf 231} (1999), 301-324.

\bibitem{ac} Cortez A., Singularités génériques des variétés de
  Schubert covexillaires, {\sl Ann. Inst. Fourier} {\bf 51} (2001),
375-393.

\bibitem{ga} Gasharov V., Sufficiency of Lakshmibai-Sandhya
  singularity conditions for Schubert varieties, {\sl Compositio Math.}
{\bf 126} (2001), 47-56.

\bibitem{klr} Kassel C., Lascoux A., Reutenauer C., The singular 
locus of a Schubert variety, prépublication IRMA, Strasbourg, mars
2001. 

\bibitem{kl} Kazhdan D., Lusztig G., Representations of Coxeter groups
  and Hecke algebras, {\sl Inventiones Math.} {\bf 53} (1979), 165-184.
 
\bibitem{ls} Lakshmibai V., Sandhya B., Criterion for smoothness
of Schubert varieties in $Sl(n)/B$, {\sl Proc. Indian Acad. Sci.}
 {\bf 100} (1990), 45-52.

\bibitem{lse} Lakshmibai V., Seshadri S., Singular locus of a
Schubert variety, {\sl Bull. A.M.S.} {\bf 11} (1984), 363-366.

\bibitem{m2} Manivel L., Le lieu singulier des vari\'et\'es de
  Schubert, preprint arXiv:math.AG/0102124.

\bibitem{man} Manivel L.,  ``Fonctions symétriques, polynômes de
    Schubert et lieux de dégénérescence'', Cours Spécialisés {\bf 3}, 
  Société Mathématique de France, 1998.


\end{thebibliography}
\end{document}